\documentstyle[twoside]{article}
\def\Vox{$^{\fbox {\/}}$}
\def\prom{
\newtheorem{prop}{Proposition}[section]
\newtheorem{cor}[prop]{Corollary}
\newtheorem{lem}[prop]{Lemma}
\newtheorem{exa}[prop]{Example}

\newtheorem{rem}[prop]{Remark}

\newtheorem{tef}[prop]{Definition}
}
\def\ieacosdo{
\author {\ \\N.\ C.\ A.\ da Costa\\F.\ A.\ Doria\\\
\\Research Group on Logic and Foundations,\\Institute for
Advanced Studies, University of S\~ao Paulo.\\Av.\ Prof.\
Luciano Gualberto, trav.\ J, 374.\\05655--010 S\~ao Paulo SP
Brazil.\\\ \\{\sc ncacosta@usp.br}\\{\sc
fadoria@rio.com.br}\\{\sc doria@lncc.br}}}
\prom

\begin {document}
\title {On the consistency of $P=NP$ with fragments of ZFC
whose own consistency strength can be measured by an ordinal
assignment.\thanks {Partially supported by FAPESP and CNPq.
Permanent address for F.\ A.\ Doria: Research Center for
Mathematical Theories of Communication and Program IDEA,
School of Communications, Federal University at Rio de
Janeiro, Av.\ Pasteur, 250. 22295--900 Rio RJ Brazil.}}
\ieacosdo
\date {May 2000\\Version 5.0}
\maketitle

\begin {abstract}
\noindent We formulate the $P<NP$ hypothesis in the case of
the satisfiability problem as a $\Pi ^0_2$ sentence, out of
which we can construct a partial recursive function $f_{\neg
A}$ so that $f_{\neg A}$ is total if and only if $P < NP$. We
then show that if $f_{\neg A}$ is total, then it isn't ${\cal
T}$--provably total (where ${\cal T}$ is a fragment of ZFC
that adequately extends PA and whose consistency is of
ordinal order). Follows that the negation of $P < NP$, that
is, $P = NP$, is consistent with those ${\cal T}$. 
\end {abstract}

\newpage

\pagestyle {myheadings}
\thispagestyle {plain}

\section {Introduction}
\markboth {da Costa, Doria}{$P = NP$}

As it is well--known, Gentzen's proof of the consistency of
arithmetic (PA) requires a transfinite induction up to
$\epsilon _0$. Ackermann's proof of the consistency of the
theory of real numbers (RT) asks for an induction up to the
first $\eta$ number, that is, the first $\epsilon$--number
$\kappa$ such that $\epsilon _{\kappa} = \kappa$. Wainer
remarked (personal communication; see also \cite {Wai,Wai1})
that some results by Kreisel \cite {Krei} can be extended to
fragments of ZFC whose own consistency is measured by an
ordinal, as in the cases of PA and RT. Therefore, it is
possible to show in this way that certain sentences are not
provable in those fragments and, as a consequence, that the
negations of such sentences are consistent with them. 

We apply those ideas to prove that $P = NP$, in the
particular case of the satisfiability problem for Boolean
expressions, is consistent with PA and other fragments of
ZFC whose consistency can be measured by an adequate 
recursive ordinal. It immediately follows that $P=NP$ in its
general form is also consistent with those theories. 

Our exposition will be somewhat informal \cite {Rog} and
intuitive; however it is easy to reformulate it in a
rigorous way. 

\subsection* {Summary of the paper}

\begin {enumerate}

\item We start from PA and from the satisfiability problem
${\cal S}$ for Boolean expressions in conjunctive normal
form (cnf).

We note $[P = NP]^{\cal S}$ the sentence that asserts that
there is a polynomial algorithm that settles all instances
of the satisfiability problem. 

\item We obtain the function $f_{\neg A}$. $f_{\neg A}$ is
total if and only if $\neg [P = NP]^{\cal S}$ (that is, $[P
< NP]^{\cal S}$) holds. 

\item We show that $f_{\neg A}$ lies ``beyond'' the Kreisel
hierarchy. Then:
\begin {enumerate}

\item If $f_{\neg A}$ is total in the standard model for PA,
then by a theorem of Kreisel \cite {Krei} PA doesn't prove
$\neg [P=NP]^{\cal S}$ and therefore $[P=NP]^{\cal S}$ is
consistent with PA. 

\item If $f_{\neg A}$ isn't total in the standard model for
PA, then again $[P = NP]^{\cal S}$ is consistent with PA, if
we suppose that PA only proves sentences that are true of
the standard model.

\item Follows that $[P = NP]^{\cal S}$ is consistent with PA.
\end {enumerate}

\item Of course this result extends to all other problems in
the $NP$--class, so that we have that $[P = NP]$ in general
is consistent with PA. 

\item The proof given for PA is also valid for any theory
that ``includes'' PA and whose consistency can be measured
by an adequate recursive ordinal. (We moreover suppose that
the theory is adequately sound.) 

\end {enumerate}

\markboth {da Costa, Doria}{$P = NP$}
\section{$[P<NP]^{\cal S}$ as a $\Pi ^0_2$ sentence}
\markboth {da Costa, Doria}{$P = NP$}

\begin {rem}\rm\label {can}
Suppose given the canonical enumeration of binary words
$$\emptyset, 0, 1, 00, 01, 10, 11, 000, 001,\ldots$$ which
code the empty word and the integers $$\sqcup \mbox {\ (\rm
blank)},  0, 1, 2,\ldots.$$ 
We take this correspondence to be fixed
for the rest of this paper. \Vox
\end {rem}

\begin {rem}\rm
We consider cnf--Boolean expressions, that is, Boolean
expressions in conjunctive normal form. 
\begin {itemize}

\item Let $x$ be a Boolean expression in cnf, adequately
coded as a binary string of length $|x|$. Let ${\sf P}_n$ be
a polynomial machine of G\"odel number $n$. (We show below
in Section \ref {ppol} how to construct a recursive sequence
of polynomial machines that suits our purposes in this
paper.) 

\item Given a string $y$ of truth--values for
the $|y|$ Boolean variables of $x$, there is a polynomial
procedure (a polynomial Turing machine which we will note
${\sf V}$) that will test whether $y$ satisfies $x$, that
is, say, ${\sf V}(\langle x,y\rangle)=1$ if and only if $y$
satisfies $x$; and equals $0$ otherwise. 

$\langle\ldots,\ldots\rangle$ is the usual \cite {Rog}
pairing function; we will only write it when required to
avoid ambiguity. 

(For the sake of completeness, we add that ${\sf V}(0,0)=1$,
that is, the empty string is satisfied by the empty string.) 

\item We formulate the recursive predicate:
$$A^*(m,x)\leftrightarrow _{\rm Def}\exists y\,({\sf
P}_m(x)=y\,\wedge\,{\sf V}(x,y)=1).$$ 
$A^*(m,x)$ is intuitively understood as ``polynomial machine
of index $m$ accepts Boolean cnf expression $x$,'' that is,
``machine $m$ inputs $x$ and outputs a line of truth values
that satisfies $x$.''

\item We can also write: 
$A^*(m,x)\leftrightarrow {\sf V}(x,{\sf P}_m(x))=1$. 

\item Form the pair $z=\langle x,y\rangle$, and let $\pi_i$,
$i=1,2$, be the usual (polynomial) projection functions.
Recall that ${\sf V}$ is a polynomial machine that inputs a
pair $\langle x,y\rangle$. Then we can consider the
predicate: 
$$\neg A(m, z)\leftrightarrow _{\rm Def} {\sf
V}(z)=1\,\wedge\, {\sf V}(\langle \pi _1 z, {\sf P}_m
(\pi _1 z)\rangle)=0,$$ 
or
$$\neg A(m,z)\leftrightarrow {\sf V}(z) = 1\,\wedge\,\neg
A^*(m,\pi_1 z).$$ 

\item $A(m,z)$ can be intuitively read as follows:
polynomial machine ${\sf P}_m$ doesn't accept the pair $z$
if and only if $z$ is such that $\pi _1 z = x$ is
satisfiable, but the output of ${\sf P}_m$ over $x = \pi_1 z$
doesn't satisfy $x$. 

This form for $\neg A(m,x)$ was suggested by F.\ Cucker to
the authors.\Vox

\end {itemize}

\end {rem}

\begin {prop}
$[P < NP]^{\cal S}\leftrightarrow\,\forall m\,\exists
z\,\neg A(m,z)$. \Vox
\end {prop}

\begin {rem}\rm
Notice that if $[P < NP]^{\cal S}$ holds, then the existence
of a single $z_0$ such that $\neg A(m, z_0)$ implies that
there are infinitely many $z'_0$ such that $\neg A(m,
z'_0)$. \Vox
\end {rem}

\subsection*{The function $f_{\neg A}$} 

\begin {tef}
$f_{\neg A}(m) =_{\rm Def} \mu _x \neg A(m,x)$. \Vox
\end {tef}

\begin {lem}
$[P<NP]^{\cal S}\leftrightarrow f_{\neg A}$ is total. \Vox
\end {lem}

\subsection*{{\rm PA}--provably total recursive functions}

The concept we now introduce originated in \cite {Krei}
Kreisel: we say that a recursive function
$f$ is {\rm PA}--provably total unary recursive if, for some
G\"odel number $e$:
\begin {enumerate}
\item {\rm PA} proves that $e$ is the G\"odel number of
$f$, and
\item For each $x$, {\rm PA} proves that the computation
of $f(x)$ converges.
\end {enumerate}

In what follows we supose that all variables are restricted
to $\omega$, the set of natural numbers. Formally,

\begin {tef} A {\rm PA}--unary function
$f:\omega\rightarrow\omega$ is {\bf {\rm PA}--provably total
unary recursive} if for some G\"odel number $e_f$ for $f$,
$${\rm PA}\vdash\forall x\,\exists
z\,(T(e_f,x,z)\wedge\,\forall y\,
(f(y)=\{e_f\}(y)).\mbox{\Vox}$$
\end {tef} 

Therefore $T(e_f,x,z)$ holds, and there is a $z$ so
that the computation of $e_f$ over $x$ stops in $z$ steps,
and not before, for every $x$. This means that we have a
proof in {\rm PA} that every computation of $f$ converges.  

\begin {rem}\rm 
Some of those non--{\rm PA}--provably total unary
functions `top' all {\rm PA}--provably total unary
functions. \Vox
\end {rem}

\begin {tef}\label {dom} For $f,g:\omega\rightarrow\omega$,
$$f \mbox{{\bf\ dominates\ }} g\leftrightarrow _{\rm Def}
\exists y\,\forall x\,(x > y\rightarrow f(x) \geq g(x)).$$ We
write $f\succ g$ for $f$ dominates $g$. \Vox
\end {tef} 

We need the next (trivial) result: 

\begin {cor}\label {over} If, for any {\rm PA}--provably
total recursive unary function $f$ we have that $g$
overshoots through $f$ infinitely many times (that is, for
infinitely many $x$,
$g(x)>f(x)$), then $g$ isn't a {\rm PA}--provably total 
unary recursive function. \Vox
\end {cor}

\markboth {da Costa, Doria}{$P=NP$}
\section{Polynomial Turing machines}\label {ppol}
\markboth {da Costa, Doria}{$P=NP$}

\begin {rem}\label {TM}\rm We describe the behavior of the
Turing machines we consider here to avoid ambiguities:
\begin {enumerate}
\item Turing machines are defined over the set $A_2^*$ of
finite words on the binary alphabet
$A_2=\{0,1\}$. 
\item Each machine has $n+1$ states $s_0, s_1,\ldots, s_n$,
where $s_0$ is the final state. (The machine stops when it
moves to $s_0$.)
\item The machine's head roams over a two--sided infinite
tape. 
\item Machines input a single binary word and either never
stop or stop when the tape has a finite, and possibly empty
set of binary words on it.  
\item The machine's {\it output word} will be the one over
which the head rests if and when $s_0$ is reached. (If the
head lies on a blank square, then we take the output word to
be the empty word.)  \Vox
\end {enumerate}
\end {rem}

\begin {rem}\rm \label {code} Our Turing machines input a
binary string $\lfloor x\rfloor$ and output a binary string
$\lfloor y\rfloor$ that stands for the numeral $y$. The
corresponding recursive functions input the numeral ${x}$
and output the numeral $y$. However as it will be always
clear from context, we write $x$ for both the binary
sequence and the numeral. \Vox
\end {rem}

\begin {rem}\rm\label {ss} We emphasize that the computable
functions we are dealing with are always given as Turing
machines. We will use upper--case sans serif
letters ({\sf M},\ldots) for the machines. If ${\sf M}_n$ is
a Turing machine of G\"odel number $n$, its input--output
relation is noted ${\sf M}_n(x)=y$.
\Vox
\end {rem}

\begin {rem}\label {0}\rm We define the {\it empty} or {\it
trivial machine} to be the Turing machine with an empty
table; we take it to be the simplest example of the identity
machine, again by definition. \Vox
\end {rem}

\subsection*{G\"odel numbering}

\begin {rem}\label{tour}\rm  Turing machines are code lines
$\xi,\xi',\ldots$, separated by blanks $\sqcup$, such as 
$\xi\sqcup\xi'\sqcup\ldots\sqcup\xi''$. Let $\Xi$ be one
such set of code lines separated by blanks. Let $\Xi '$ be
obtained out of $\Xi$ by a permutation of the lines $\xi,
\xi ',\ldots$. Both $\Xi$ and $\Xi '$ are seen as different
machines that compute the same algoritmic function, that is,
in this case, if $f_{\Xi}$ ($f_{\Xi '}$) is computed by $\Xi$
($\Xi '$), then $f_{\Xi} = f_{\Xi '}$. \Vox
\end {rem}

We allow some freedom in the construction of those `code
lines' $\xi\ldots$. Usual program instructions are
acceptable. 

\begin {rem}\rm 
\label {Ind} We can recursively enumerate all tables for
Turing machines as described in Remarks \ref {TM} and \ref
{tour}. We suppose that our G\"odel numbering for Turing
machines arises out of the following: 
\begin {itemize}
\item We list all words in the alphabet used to describe the
Turing machine tables. 
\item The set of tables is recursive, and so, given a binary
word $x$,
\begin {enumerate}
\item If it translates into a Turing machine table through
the usual 1--1 correspondence between words and machines,
then the corresponding numeral codes that machine. 
\item If not, we impose that it will code the trivial
machine. \Vox
\end {enumerate} 
\end {itemize}
\end {rem}

\begin {rem}\label {26}\rm 
We can explicitly construct Turing machines ${\sf
F}_{\alpha}$, $0\leq\alpha < \epsilon _0$ for the
Kreisel hierarchy of fast--growing total recursive
functions $f_{\alpha}$, which are used here according to the
prescriptions in \cite {Krei,Spen,Wai,Wai1}. \Vox 
\end {rem} 

We present our arguments and constructions for PA, that is,
for the machines ${\sf F}_0, {\sf F}_1, \ldots, {\sf
F}_{\alpha},\ldots$,
$\alpha < \epsilon _0$. Extension to other theories whose
consistency is of ordinal rank will be discussed in the
course of our argument. 

\subsection*{Machines bounded by a clock}

In order to deal with the set of all polynomial
machines we must resort to a trick. That is, we deal with a
recursive set of {\it expressions} for polynomial machines,
so that for each ``concrete'' polynomial machine there is an
expression for that machine, and no expression in the
listing represents a nonpolynomial machine. 

\begin {rem}\rm  The idea goes as follows: 
\begin {enumerate}
\item We use at first a variant of the \cite {Bak}
Baker--Gill--Solovay trick: we consider all couples $\langle
{\sf M}_n, {\sf C}_p\rangle$ where ${\sf M}_n$ is a Turing
machine, and ${\sf C}_p$ is a polynomial clock (Definition
\ref {Kl}) that, for an input $x$ of length
$|x|$, shuts down the machine after ${\sf C}_p(|x|)$ steps,
if the machine is still running \cite {Bak}. 

We are soon going to specify the structure of the clocks we
use; they will always be PA--provably total Turing machines. 
\item Since there is an uniform recursive procedure to obtain
a Turing machine out of each pair $\langle {\sf M}_n, {\sf
C}_p\rangle$, we `embed' this recursive sequence of
pairs (that represent polynomial machines) into the sequence
of all Turing machines. 
\item Notice that each polynomial machine will be
`represented' by several Turing machines.
\item We will take the Turing machine G\"odel number as the
G\"odel number for our polynomial machines. 
\end {enumerate}

In order to handle the machine pairs, we define below an
adequate indexing system for them that includes both the
machine and the clock that bounds its operation time. \Vox
\end {rem}

\subsection*{Parametrized polynomial clocks}

\begin {rem}\rm We suppose that the clock acts
as follows. If ${\sf M}_n(x)$  stops before the clock
determines it to shut down, the output is precisely ${\sf
M}_n(x)$. However if the bound in the number of processing
steps is reached before ${\sf M}_n(x)$ stops, the clock stops
it and orders its state to move to $s_0$. The
machine's output is then agreed to be $0$. \Vox
\end {rem}

\begin {tef}\label {Kl} A {\bf parametrized polynomial
clock} ${\sf C}_{{\sf F}(k)}$ is a total Turing machine
that depends on a positive integer $k$ and on a total
recursive ${\sf F}$, and that satisfies:
\begin {enumerate} 

\item The clock is coupled to another Turing machine
${\sf M}_n$. 

\item Whenever $x$ is input to ${\sf M}_n$, $x$ is also
input to ${\sf C}_{{\sf F}(k)}$. 

\item ${\sf C}_{{\sf F}(k)}$ computes $|x|^{{\sf F}(k)} +
{\sf F}(k)$. 

\item If ${\sf M}_n(x)$ hasn't stopped before, then ${\sf
C}_{{\sf F}(k)}$ shuts it down after it has completed
$|x|^{{\sf F}(k)} + {\sf F}(k)$ steps in the computation
over $x$ and makes it output $0$. \Vox 
\end {enumerate}
\end {tef}

\begin {rem}\rm
We will restrict our attention to the following class of
clocks that will be coupled to the Turing machines in order
to make them polynomial machines:
\begin {enumerate}
\item ${\sf C}_p$ is the clock that shuts down the Turing
machine ${\sf M}_n$ in the machine pair $\langle {\sf
M}_n,{\sf C}_p\rangle$ after $|x|^p + p$ steps, as described
above.
\item \label {cu} Now consider the Kreisel hierarchy
\cite {Sol,Krei,Wai}. For the ordinal indices $$\omega,
\omega^{\omega},\omega^{\omega ^{\omega}},\ldots <\epsilon
_0,$$ in its obvious enumeration, form the corresponding
functions ${\sf F}_0, {\sf F}_1, \ldots$, where `$0$' stands
for $\omega$, `$1$' stands for $\omega ^{\omega}$,\ldots,
and so on. 

\item Then ${\sf C}_{{\sf F}_{\alpha}(k)}$ denotes the clock
that shuts down the Turing machine in the manner prescribed
above after $x^{{\sf F}_{\alpha}(k)} + {{\sf F}_{\alpha}(k)}$
steps, where ${\sf F}_{\alpha}$ is indexed as described. 

\item We will consider two classes of clocks: those that
correspond to the polynomials $|x|^p + p$, $p = 0, 1,
2,\ldots$, and the (families of) clocks that correspond to
the parametrized polynomials $x^{{\sf F}_{\alpha}(k)} +
{{\sf F}_{\alpha}(k)}$. 

\item We suppose given a fixed recursive enumeration for
those two kinds of clocks. 
\end {enumerate}

This restriction won't affect our result. \Vox
\end {rem}

In what follows we will only deal with those two kinds of
clocks. 

\begin {rem}\rm
For the case where we are dealing with a theory $\cal T$ of
consistency rank $\zeta > \epsilon _0$, recall that $\zeta$
is a constructive ordinal. Therefore we can always find a
recursive subset of ordinals as those in step \ref {cu}
above to proceed as stated. \Vox
\end {rem}

\begin {rem}\rm Again let 
$\langle\ldots ,\ldots\rangle$ denote the usual
\cite {Rog} degree--2 polynomial recursive 1--1 pairing
function
$\langle\ldots,\ldots\rangle:\omega\times\omega\rightarrow\omega$.
\Vox
\end {rem}

We define: 

\begin {tef}\label {X} ${\sf
P}_{p}=\langle {\sf M}_i,{\sf C}_j\rangle$, ${\sf C}_j$ a
parametrized clock. Call their set ${\cal P}$, the set of
Turing machines coupled to parametrized polynomial clocks.
\Vox
\end {tef}

\begin {rem}\label {Tt} \rm 
${\cal P}$ contains all possible ordered pairs as in the
preceding definitions; the number $p$ can be intuitively
seen as coding an expression for a Turing machine with a
clock. In what follows ${\cal M}$ will denote the sequence
of all Turing machines. \Vox
\end {rem} 

\begin {prop}\label {rek} There is a recursive 1--1
embedding $\sigma :{\cal P}\rightarrow {\cal M}$ of the
polynomial machines represented by pairs $\langle {\sf
M}_m,{\sf C}_n\rangle$, into polynomially--bounded Turing
machines given by their tables. The set $\sigma {\cal P}$ is recursive in 
$\cal
M$. \Vox
\end {prop}

\begin {rem}\label {reason}\rm 
Suppose that we are given the image $\sigma{\cal P}$ of that
map into ${\cal M}$, and suppose moreover that we add to
that image all explicitly defined polynomial machines that
we use in this construction, such as the ones of the form
given in Lemma \ref {mak} below. Those machines are
finite--output machines (in a sense made clear in that
lemma); we note their set ${\cal F}$. 

The resulting set ${\cal F}\cup\sigma {\cal P}$ of polynomial
machines is recursive, and the consistency result we
describe here related to $f_{\neg A}$ is valid if and only if
it is valid when extended both to that subset $\sigma {\cal
P}\subset{\cal M}$ when coded by their own Turing--machine
G\"odel numbers and to ${\cal F}\cup\sigma {\cal P}$. \Vox
\end {rem}

\begin {rem}\rm\label {reason1} 
We form an extended $f'$ by defining it to be equal to
$f$ over $\sigma {\cal P}$ and 
$0$ otherwise. Therefore $f_{\neg A}$ is total if and only
if the extension
$f'_{\neg A}$ is total. \Vox
\end {rem}

Since there is such a recursive, uniform procedure $\sigma$
that allows us to obtain actual tables for polynomial Turing
machines given an arbitrary machine and clock as described
here, we can form their compositions, which are again
polynomial machines, granted a recursive rule to define a
bounding clock. 

We can easily obtain, out of an adequate definition for the
composition operation, that:

\begin {prop}
${\cal P}$ is closed under composition, that is, if ${\sf
P}, {\sf P}'\in {\cal P}$, then we can obtain a ${\sf
P}''\in {\cal P}$ such that ${\sf P}''={\sf P}\circ{\sf
P}'$ as a function. \Vox
\end {prop}

\subsection*{Two side comments}

(See \cite {MY}.) Polynomial Turing machines can be very
powerful:

\begin {rem}\rm\label {Exa1}

We first show that if ${\sf G}$ is a fast--growing
(``hard,'' ``intractable'') recursive function that is the
most efficient way of settling a given problem, and if ${\sf
H}$ is another such function, then there is a family of
polynomial Turing machines parametrized by $n\in\omega$ such
that machine $n$ settles exactly ${\sf H}(n)$ instances of
the problem we are dealing with: 

\begin {lem}\label {mak} Let $n$ be any fixed natural number
and let
${\sf G}$, ${\sf H}$ be Turing machines that compute two
arbitrary, monotonic increasing, total unary recursive
functions. Then for ${\sf P}_{({\sf G},{\sf H})}$ given by:
\begin {enumerate}
\item ${\sf P}_{({\sf G},{\sf H})}(x) = {\sf G}(x)$, $x\leq
{\sf H}(n)$. 

\item ${\sf P}_{({\sf G},{\sf H})}(x) = 0$, $x>{\sf H}(n)$. 
\end {enumerate} there is a polynomial algorithm for it. 
\end {lem}

{\it Sketch of proof}\,: We can construct the table of ${\sf
P}$ in such a way that the operation time of ${\sf P}$ is
bounded by a constant larger than the amount of time
required for the largest computation of ${\sf G}(x)$, $x\leq
{\sf H}(n)$. \Vox

\

\noindent (We suppose that $0$ is no solution for the 
problem which ${\sf G}$ settles, for $x>0$.)  So, we have
established the existence of the family that we are looking
for. \Vox
\end {rem}

An interesting intractable problem for polynomial machines is
described below:

\begin {lem} It is ``hard'' to determine, for an arbitrary
$n$, the value of ${\sf H}(n)$, which is the cardinality of
the set of instances of the desired problem which are
settled by the $n$--th machine in that family. \Vox
\end {lem}

\begin {tef}
The set of all machines as in Lemma \ref {mak} is noted
${\cal F}$. \Vox
\end {tef}

\markboth {da Costa, Doria}{$P = NP$}
\section{Sketch of the main result}
\markboth {da Costa, Doria}{$P = NP$}

\begin {rem}\rm 
>From here on we act according to the following: 
\begin {itemize}
\item  We have added to the (recursive) image $\sigma{\cal
P}\subset {\cal M}$ all machines in Remark \ref {Exa1}
besides the images of the pairs $\langle {\sf M}, {\sf
C}\rangle$ to get ${\cal F}\cup\sigma {\cal P}$. 

\item That $f_{\neg A}$ has been extended 
to $f'_{\neg A}$ over ${\cal M}$ as in Remarks \ref {reason}
and \ref {reason1}. However by an abuse of language we
will use $f_{\neg A}$ for the extended function. \Vox
\end {itemize}
\end {rem}

We can now summarize our main argument. 

\begin {rem}\rm 
Recall that if $f$ is PA--provably total recursive, then
there is a PA--provably total recursive $g$ such that
$g\succ f$. 

$f_{\neg A}$ will be seen not to be PA--provably total. \Vox
\end {rem}

\begin {rem}\rm\label {MOT} We add more detail to our
argument. Suppose that we have constructed $\sigma {\cal
P}\subset {\cal M}$, and that $f_{\neg A}$ has been extended
as indicated above to the whole of ${\cal M}$: 
\begin {itemize}
\item Let ${\sf F}$ be a Turing machine that computes a 
recursive, unary total function which dominates all
PA--provably total functions. We can take ${\sf F}= {\sf
F}_{\epsilon _0}$ in one of the hierarchies described in
\cite {Sol,Krei,Spen,Wai}. 

\item Let ${\sf E}$ denote a fixed, explicitly given
exponential Turing machine that solves any instance of the
satisfiability problem. 

\item For each natural number $n$, form a polynomial machine
${\sf Q}^{{\sf F}(n)}$ that operates as follows: 
\begin {enumerate}
\item Put $K= {\sf F}(n)$. 

\item For $x \leq K$,
${\sf Q}^{{\sf F}(n)}(x)={\sf E}(x)$. 

\item For $x > K$, ${\sf Q}^{{\sf F}(n)}(x) = 0$. 
\end {enumerate}

(Notice that $0$ is never a solution, for $x > 0$.) So up to
${\sf F}(n)$ the machine ${\sf Q}^{{\sf F}(n)}$ equals ${\sf
E}(x)$. From then on it always prints $0$. 

This machine ${\sf Q}^{{\sf F}(n)}$ is polynomial (see
Lemma \ref {mak}). 

\item Its Turing machine table can be described out of the
following set of instructions:
\begin {enumerate}

\item The code for the parameter $n$. (Perhaps a single
instruction line, as $y = n$.) 

\item The code for computing ${\sf F}$. (Result of the
computation should be ${\sf F}(y)$, for $y$ as input.) 

\item Instructions for the computation of ${\sf Q}^z(x)$, for
$z = {\sf F}(y)$. (Here the instructions involve
the variables $x$ and $y$.) This involves the fixed code
for ${\sf E}$. 
\end {enumerate}

(See in Remark \ref {tour} the coding procedures for Turing
machines that we use in this paper.) 

\item Notice that (if we use the coding techniques described 
in Remark \ref{tour}), for each $n$, there are constant
natural numbers $a, b$ such that the G\"odel number for
${\sf Q}^{{\sf F}(n)}$ equals $an + b$. (See Proposition
\ref {shoo}.) 

\item The operation time of ${\sf Q}^{{\sf F}(n)}$ is bounded
by a polynomial clock with bound ${\sf F}'(n) + x^{{\sf
F}'(|n|)}$, each
$n$, where $|x|$ is the length of the binary word $x$, for
sufficiently large ${\sf F}'>{\sf F}$. (It is enough to take
${\sf F}' = {\sf F}_{\epsilon _0 + 1}$.) 

\item Again the instructions for the operation of that clock
are such that, given each $n$, their G\"odel numbers are
given by $a'n + b'$. (See Proposition \ref {shoo}.) 

\item The pairing $\langle an + b, a'n + b'\rangle$ is 
quadratic on $n$ \cite {Rog}. 

\item Thus (intuitively) $f_{\neg A}(\langle an + b, a'n +
b'\rangle)= {\sf F}(n)+ 1$ (over the machine pairs), or
$f_{\neg A}(an + b) = {\sf F}(n) + 1$ over ${\cal
F}\cup\sigma {\cal P}$. Intuitively again, our function
$f_{\neg A}$ overshoots infinitely many times through every
provably total unary recursive function in PA. 

Given the preceding results, plus Lemma \ref {DDom}, we
conclude that $f_{\neg A}$ isn't dominated by any such
PA--provable total recursive function. 
\end {itemize}

See Remark \ref {phoo} in the final Section of the paper.
\Vox
\end {rem} 

\begin {rem}\rm 
We can also argue as follows: instead of using the ``ceiling
function'' ${\sf F}_{\epsilon _0}$, let us be given the
Kreisel \cite {Sol,Krei} hierarchy $\{{\sf F}_0, {\sf
F}_1,\ldots,{\sf F}_{\beta},\ldots\}$, $0\leq\alpha
<\epsilon _0$, of PA--provably total functions and let 
${\sf F}_{\alpha}$ be in that hierarchy. Then:
\begin {itemize}
\item  We make
the preceding construction for a ${\sf F}_{\beta}$ such that
${\sf F}_{\beta}\succ {\sf F}_{\alpha}$.
\item This will be true even if ${\sf F}_{\alpha}$ is
composed with a quadratic function of its variable, that is,
${\sf F}_{\beta}\succ {\sf F}_{\alpha}\circ
u$, where $u$ is quadratic, and so on. 
\item This can be repeated for any $\alpha$ in the Kreisel
hierarchy. 
\item So, we conclude that no PA--provably total recursive
function can dominate $f_{\neg A}$. \Vox
\end {itemize}
\end {rem} 

\begin {rem}\rm
Since we will work within the set of Turing machines ${\cal
M}$, as given in Proposition \ref {rek} the G\"odel number
for the machines ${\sf Q}$ will be a linear function
$an + b$. \Vox
\end {rem} 

\begin {rem}\label {reason2}\rm 
We do not obtain this result if we restrict ourselves to an
enumeration such as the one in \cite {Bak}. (We thank S.\
Wainer for that observation.) 

Yet we feel it isn't natural to exclude polynomial Turing
machines such as those in our family ${\sf Q}^{{\sf F}(n)}$,
any $n$ (Example \ref {Exa1}) since that family contains
rather simple polynomial machines with the property that it
is very hard to compute the input values at which they start
to output zeros forever.

The gist of the matter is this fact: a set of machines (given
by their tables) can only contain polynomial Turing machines,
but it may be very, very hard to decide some of its
properties. Our family ${\sf Q}^{{\sf F}(n)}$, all $n$, is
an example of such a set.
\Vox
\end {rem}

\markboth {da Costa, Doria}{$P = NP$}
\section{Computation of some infinite segments of the
function $f_{\neg A}(n)$}
\markboth {da Costa, Doria}{$P = NP$}

We can restrict our attention to the machines in ${\cal
F}$. We will argue for PA, but it is easy to extend our
results to the ${\cal T}$ already characterized. We recall:

\begin {tef}\label{GG} 
 $f_{\neg A}(m)=_{\rm Def}\mu _x \neg A(m,x)$. \Vox
\end {tef} 

\begin {rem}\rm\label {hypo}
Now either ``$f_{\neg A}$ is total'' is true of the standard
model or it isn't. \Vox
\end {rem}

We prove here:

\begin {prop}\label {main}
Given any ${\sf F}_{\alpha}$, $0\leq\alpha <\epsilon _0$
then for no $\alpha$ does ${\sf F}_{\alpha}$ dominate
$f_{\neg A}$. 
\end {prop}

\begin {cor}
If {\rm PA} is arithmetically sound, then ${\rm PA}\not\vdash
f_{\neg A}$ {\em is total.} Therefore $[P = NP]^{\cal S}$ is
consistent with {\rm PA}. \Vox
\end {cor} 

{\it Proof}\,: Will be given throughout this Section. 

\begin {rem}\rm
For ``arithmetic soundness'' (or ``arithmetic consistency'')
see \cite {CosDo1}. \Vox
\end {rem}

\subsection*{Computation of the G\"odel number of ${\sf
Q}^{{\sf F}_{\alpha}(n)}$} 

\begin {rem}\rm
We must clearly separate the two kinds of machine
codes that we will be using for the purposes of our proof:
\begin {itemize}
\item The G\"odel number of a Turing machine is its ${\cal
M}$--code, index or simply its G\"odel number. 
\item We also code our polynomial machines out
of the G\"odel numbers of the paired machine--and--clock
$\langle {\sf M}, {\sf C}\rangle$ representation. This
specific code for polynomial machines is their ${\cal
P}$--index or code. \Vox
\end {itemize}
\end {rem}

\begin {rem}\rm
The argument presented in this subsection is very simple. 
G\"odel numbers for the machines in our proof will turn out
to be as in an arithmetic progression like  
$547$, $647$, $747$, $847,\ldots$. 

To add some more detail: 
\begin {itemize}
\item {\it G\"odel numbers.} We show that the G\"odel
numbers for the machines that interest us and the
corresponding clocks can be written as a string 
$[n] s$ ($n = 0, 1, 2,\ldots$), $[n]$ is a binary
string that includes the binary string for $n$, and $s$ is a
fixed binary string; concatenation of both strings is
indicated by their juxtaposition. As we noticed, those
G\"odel numbers are in an arithmetic progression.

\item {\it Indices for machine pairs.} As the usual pairing
function \cite {Rog}
$\langle x,y\rangle$ is a degree--2 polynomial, the code for
a pair $\langle {\sf M}, {\sf C}\rangle$ will depend on a
degree--2 polynomial on $n$ when parametrized as above. \Vox
\end {itemize}
\end {rem}

Our goal is to compute the ${\cal M}$-- and ${\cal
P}$--indices of the machines  ${\sf Q}^{{\sf
F}_{\alpha}(n)}$ as a function of the tables for ${\sf
F}_{\alpha}$, for those finite machines described above, and
as a function of $n$. 

For notation see Remarks \ref {tour} and \ref {Ind}. We
consider machines ${\sf Q}^{{\sf F}_j(n)}$ and the pairs
$\langle {\sf Q}^{{\sf F}_j(n)}, {\sf C}\rangle$, with a
clock ${\sf C}$ that bounds the polynomial machine without
cutting it off before it stops by itself. 

We will first deal with the indexing over ${\cal M}$, and
then over ${\cal P}$. See Remark \ref {Ind} and
recall that ``garbage'' is mapped on the trivial machine. 

\begin {rem}\rm
For the remainder of this Section, Turing machine tables
will be given in the form of sketchy, summary programs. \Vox
\end {rem}

\subsection*{Estimates for G\"odel numbers}

Our result on G\"odel numbers is: 

\begin {prop}\label {shoo} We can write the table for ${\sf
Q}^{{\sf F}_{\alpha}(n)}$, each $n$, so that its G\"odel
number (${\cal M}$--index) $\rho (n,\alpha) = a_{\alpha} +
(2^{q_{\alpha}}-1)n$, $a_{\alpha}, q_{\alpha}$ positive
integers, $q_{\alpha}$ a constant that depends on the
G\"odel number for ${\sf F}_{\alpha}$. 
\end {prop} 

\begin {exa}\rm
We can give an example to make things explicit. The (sketchy)
program we use will look like: 
\begin {enumerate}
\item\label {1'} Start. 
\item\label {2'} $y = n$.
\item\label {3'} Input program for ${\sf F}_{\alpha}$.

Compute ${\sf F}_{\alpha}(y)$. 
\item\label {4'} Input program for $2^x$. 

Compute $2^z$, for $z = {\sf F}_{\alpha}(y)$. 
\item\label {5'} Output $2^z+1$. 
\item Stop. 
\end {enumerate}
\end {exa}

{\it Proof}\,: Recall that $\lfloor x\rfloor$ is a binary
string for $x$. We may write the program for the polynomial
Turing machine 
${\sf Q}^{{\sf F}_{\alpha}(n)}$ as the concatenation
$\lfloor [n]\rfloor \lfloor \xi
'_{\alpha}\rfloor\lfloor\xi''\rfloor$, where:
\begin {itemize}
\item  $[n]$ codes lines \ref {1'} and \ref {2'}.

It is a binary string that includes the bits for the numeral
$n$ (that is the only variable portion in the machine's
program). 
\item $\xi'_{\alpha}$ codes line \ref {3'}; essentially the
program for ${\sf F}_{\alpha}$;
\item $\xi''$ codes lines \ref {4'} and \ref {5'}. 
\end {itemize}

Given our G\"odel numbering conventions, the value of the
G\"odel number for ${\sf Q}^{{\sf F}_{\alpha}(n)}$,
for $q_{\alpha}=\lfloor\xi'_{\alpha}\xi''\rfloor$, is
given by a binary number $\rho (n,\alpha) = \lfloor [n]
q_{\alpha}\rfloor$. 

Machines ${\sf Q}^{{\sf F}_{\alpha}(n)}$, for each $n$, as
given by the tables so described are regularly spaced
among the ${\sf M}_m$ in the ordering we have given for
${\cal M}$ (see Remark \ref {Ind}), that is, at most 
$2^{q_{\alpha}}-1$ machines lie between machine
${\sf Q}^{{\sf F}_{\alpha}(n-1)}$ and machine ${\sf Q}^{{\sf
F}_{\alpha}(n)}$ in that arrangement. In other words, the
G\"odel number of a ${\sf Q}^{{\sf
F}_{\alpha}(n)}$ machine,
$\rho (n,\alpha) = a_{\alpha}+(2^{q_{\alpha}}-1)n$,
constants $a_{\alpha},q_{\alpha}$ positive integers. \Vox 

\begin {rem}\rm Now we must estimate the G\"odel number of a
(reasonably small) clock ${\sf C}_{{\sf F}_p}$ that bounds
${\sf Q}^{{\sf F}_{\alpha}(n)}$. Out of the preceding
argument we see that the table for the clock has the form:
$$[n][\mbox{\ program for ${\sf
F}_{\alpha}$\ }][\mbox{\ coupling instructions.}]$$ 

The only variable portion of it is $n$. Again we have, for
the G\"odel numbers $K(\alpha,n)$ of the clocks that bound
the ${\sf Q}^{{\sf F}_{\alpha}(n)}$, that, for each
$n$, that they equal $b_{\alpha} + (2^{q'_{\alpha}}-1)n$,
$q'_{\alpha}, b_{\alpha}$, positive constants. \Vox
\end {rem}

\begin {rem}\rm Then the ${\cal
P}$--index for the couple
$$\langle {\sf Q}^{{\sf F}_{\alpha}(n)},
{\sf C}_{{\sf F}_p}\rangle$$ is a degree--2 polynomial on
those linear functions. That is, the series of G\"odel
numbers for those machines in the representation (set of
clock--bounded machines) ${\cal P}$ is bounded by a very
reasonable function. 

However the elements of $\sigma {\cal P}\subset {\cal M}$ are
linearly spaced in ${\cal M}$ when this latter set is
ordered by the machines' G\"odel numbers. For in order to
make $\sigma {\cal P}$ recursive we write the couples
$\langle {\sf Q}^{{\sf F}_{\alpha}(n)}, {\sf C}_{{\sf
F}_p}\rangle$ as [parameter for the clock][instructions for
the clock][instructions for the machine]. The first set of
strings may be empty or is easily recognizable; the next one
is also a recognizable set (the clock); and then we have an
arbitrary machine. The first two sets provide the
identification for the members of $\sigma {\cal P}$ in
${\cal M}$. 

Then the argument presented above shows that those machines
are linearly spaced. \Vox 
\end {rem}

\begin {cor} The ${\cal P}$--index $\lambda
\langle {\sf Q}^{{\sf F}_{\alpha}(n)},
{\sf C}_{{\sf F}_p}\rangle$ is a degree--2 polynomial on $n$.
\Vox
\end {cor} 

\begin {cor} The G\"odel number $\sigma {\sf Q}^{{\sf
F}_{\alpha}(n)}$ is linear on $n$. \Vox
\end {cor} 

\subsection*{Fast--growing functions ``embedded'' into
$f_{\neg A}$}

We obtain a kind of ``copy'' of a fast--growing function
${\sf F}$ within $f_{\neg A}$ in such a way that the the
images
$\sigma {\sf F}(n)$, $\sigma {\sf F}(n+1),\ldots$, are
separated in a ``controlled'' way within
$f_{\neg A}$, that is, ${\sf F}(n) = f_{\neg A}(p(n))$, ${\sf
F}(n+1) = f_{\neg A}(p(n+1)),\ldots$, where $p$ is a
polynomial.   

\subsection*{The main lemma}

\begin {lem}
For no ${\sf F}_{\alpha}$ does ${\sf F}_{\alpha}$ dominate
$f_{\neg A}$.
\end {lem}

\subsection*{Proof of the main lemma}

Let the ${\sf F}_{\alpha}$, $\alpha$ a countable ordinal,
$\alpha < \epsilon _0$, be the dominating functions in the
Kreisel hierarchy \cite {Sol,Krei,Spen,Wai} in PA. Now recall
(Remarks \ref {tour} and \ref {Ind}) that, depending on the
unicity or not of $\rho {\sf X}$, the index (G\"odel number)
of machine ${\sf X}$ in PA, it is uniquely defined but for
the trivial machine. (See Definition \ref {tour}.) Then:

\begin {tef} For any $\alpha$ define the map from 
${\sf F}_{\alpha}$ to $f_{\neg A}$ given by:
$${\sf F}_{\alpha}(n)\mapsto f_{\neg A}(\sigma {\sf
Q}^{{\sf F}_{\alpha}(n)}) = {\sf F}_{\alpha}(n) + 1.\mbox{\
\Vox}$$ 
\end {tef}

(We don't need the clock here, as we have added the machines
described in Example \ref {Exa1} to $\sigma {\cal P}$.) 

Put $\psi_{\alpha} (n)=\sigma {\sf Q}^{{\sf
F}_{\alpha}(n)}$; it is at worst linear on $n$. We will be
interested in the ``peaks'' of
$f_{\neg A}(m)$ at
$m=\psi (n)$, $n = 0, 1, 2, \ldots$. 

Restricted to those values of the variable $n$, $$f_{\neg
A}(\psi_{\alpha} (n)) = {\sf F}_{\alpha}(n) + 1.$$ 

For all $\beta <\alpha$, $f_{\neg A}(n)$ overshoots through
${\sf F}_{\beta} (n)$ infinitely many times. Let's see how
it is done. We first need:

\begin {lem}\label {DDom}
For any positive--definite
polynomial $p$ and any $\alpha < \epsilon _0$, there is a
$\beta$, $\alpha <\beta < \epsilon _0$ such that
${\sf F}_{\beta}\succ {\sf F}_{\alpha}\circ p$.
\end {lem}

{\it Proof}\,: We use this simple characterization for the
Kreisel hierarchy \cite {Spen}:  
\begin {enumerate}
\item ${\sf F}_0 = 0$, ${\sf F}_1 = 2x$. 

\item ${\sf F}_{\beta + 1}(x) = {\sf F}_{\beta}^{(x)}(1)$.

\item ${\sf F}_{\beta}(x) = {\sf F}_{\beta (x)}(x)$, where
the sequence of ordinals \ldots $\beta (x)$\ldots converges
to the limit ordinal $\beta$. 
\end {enumerate} 

Now suppose that $p$ is a $(n - 2)$--degree positive definite
polynomial. Then: 

$${\sf F}_{\beta + n + 1}(x) = {\sf
F}_{\beta + 1}^{(x^n)}(1)\succ {\sf F}_{\beta}^{(p(x))}(1).$$

As $x^n\succ p(x)$, we are done. For the limit ordinal case,
use induction over the sequence of ordinals that converge to
limit $\beta$. \Vox

\

Now to conclude our proof: 

\begin {itemize}
\item Suppose that ${\sf F}_{\beta}\succ f_R$. 

\item From Lemma \ref {DDom} there is an $\alpha$ such that
${\sf F}_{\alpha}\succ {\sf F}_{\beta}\circ p$, $p$ a fixed 
polynomial. 

\item More precisely, for $x> x_0$, ${\sf F}_{\alpha}(x) >
{\sf F}_{\beta}(p(x))$, or, given a change of variables,
${\sf F}_{\alpha}(p^{-1}(y)) > {\sf F}_{\beta}(y)$. 

\item Finally, from our construction, we have that
$f_{\neg A}(n) = {\sf F}_{\alpha}(p^{-1}(n)) + 1 > {\sf
F}_{\beta} (n)$, $n > n_0$, and thus our proof. \Vox
\end {itemize} 

This is of course only valid for $f_{\neg A}$ restricted to
the previously given values of $\psi_{\alpha} (n)$; we aren't
interested here in what happens to the in--between values. 
So we conclude: 

\begin {lem} For no $\alpha < \epsilon _0$ does ${\sf
F}_{\alpha}$ dominate $f_{\neg A}$.  
\end {lem}

{\it Proof}\,: By construction and from Proposition \ref
{shoo}. \Vox

\begin {rem}\rm  As already spelled out in the Remark that
opens this Section, the idea is that $f_{\neg A}$ contains
``peaks'' that are spaced in a controlled, bounded way.
Those peaks overtake infinitely many times the monotonic
functions ${\sf F}_{\alpha}$. \Vox
\end {rem}

\subsection*{The consistency result}

Recall:

\begin {tef}
A predicate $P(x_0,\ldots,x_n)$ is {\rm PA}--{\bf provably
recursive} if its characteristic function is
{\rm PA}--provably recursive. \Vox
\end {tef}

\begin {cor}
$A(m, x)$ is {\rm PA}--provably recursive. \Vox
\end {cor}

We also need \cite {Krei}:

\begin {prop}\label {KKKo} 
If $P(x, y)$ is {\rm PA}--provably recursive, given $$f_P=\mu
_y P(x,y),$$ if ${\rm PA}\vdash\forall x\,\exists y\,
P(x,y)$, then there is an $\alpha<\epsilon _0$ such that
${\sf F}_{\alpha}\succ f_P$. \Vox
\end {prop} 

For adequately sound PA, and given Remark \ref {hypo}: 

\begin {cor}
${\rm PA}\not\vdash\forall x\,\exists y\, \neg A(m,x)$. \Vox
\end {cor}

\begin {cor}\label {kons}
${\rm PA}\not\vdash\neg [P = NP]^{\cal S}$. \Vox
\end {cor}

\begin {rem}\rm
The soundness condition we require was called ``arithmetic
consistency,'' or ``arithmetic soundness,'' in some previous
papers by the authors. See, for instance, \cite {CosDo1}.
\Vox
\end {rem}

\begin {rem}\rm 
Recall that, as already indicated, the result in Corollary
\ref {kons} extends to any fragment ${\cal T}$ of {\rm ZFC}
that includes {\rm PA} and whose consistency strength is
measured by a countable, recursive, ordinal. \Vox
\end {rem} 

\subsection*{$f_{\neg A}$ is self--similar}

\begin {rem}\label {ssim}\rm
$f_{\neg A}$ is self--similar in the following sense:
intuitively, we can embed $h$ into itself, as in the
previous proof. Thus there will be countably infinite many
copies of $f_{\neg A}$ within itself. The same is true of any
recursive function in whose construction $f_{\neg A}$
appears. \Vox
\end {rem}

\subsection*{Extension to theories beyond {\rm PA} which are
not of ordinal consistency rank}

\begin {rem}\rm\label {phoo}
We are interested in the further extension of those results
to the whole of ZFC. Within ZFC there are functions that are
total recursive but not ZFC--provably total recursive
(granted that one  supposes that every arithmetic formula
that is provable in ZFC is true of the standard model).
However it isn't immediately clear how we can extend our
preceding result to the whole of ZFC, as the Kreisel
hierarchy \cite {Krei,Wai} stops at $f_{\epsilon _0}$. 

Proposition \ref {KKKo} applies to theories for which an
``ordinal assignment''  measure of consistency--strength is
known. Although it applies to many stronger theories than
PA, it cannot apply  to ZF or ZFC because we have no idea
how to measure its consistency  strength. 

Any use of ZF or ZFC in this context is highly confusing
because  of this lack of ``ordinal measures.'' (Personal
communication by S.\ S.\ Wainer.) \Vox
\end {rem}

\markboth {da Costa, Doria}{$P = NP$}
\section{Acknowledgments}
\markboth {da Costa, Doria}{$P = NP$}

We must thank S.\ S.\ Wainer who very kindly communicated to
us some of his recent work about the hierarchy of total
recursive functions; M.\ Benda also pointed out to the
authors some recent work on similar fast--growing
functions.  We also thank F.\ Cucker \cite {Cucker} for a
suggestion about $A(m,x)$ and for an exchange of e--mails on
a variant of $f_{\neg A}$. 

The authors wish to thank  FAPESP, CNPq and CAPES, in their
respective Philosophy Sections, as well as the Research Group
on Logic and Foundations at the IEA/USP (as well as the
supporting team at the IEA), for grants related to the
present work. Office space was given to FAD by Program
IDEA and its chairman Marcio Tavares d'Amaral, whom we
gratefully thank. Partial support for FAD and computer
formatting of this paper are due to the Project {\it Casa da
Torre}. 

\markboth {da Costa, Doria}{$P = NP$}

\bibliographystyle {plain}

\begin {thebibliography}{99}
\bibitem {Bak} T.\ Baker, J.\ Gill, R.\ Solovay,
``Relativizations of the $P=?NP$ question,'' {\it SIAM
Journal of Computing} {\bf 4}, 431 (1975). 
\bibitem {Cucker} F.\ Cucker, e--mail messages to the
authors (1999). 
\bibitem {CosDo1} N.\ C.\ A.\ da Costa and F.\ A.\ Doria,
``Undecidability and incompleteness in classical
mechanics,'' {\it International Journal of Theoretical
Physics} {\bf 30}, 1041 (1991). 
\bibitem {CosDo2} N.\ C.\ A.\ da Costa and F.\ A.\ Doria,
``On total recursive but not ${\cal T}$--total recursive
functions,'' preprint IEA--USP (1999). 
\bibitem {Do} F.\ A.\ Doria, ``Is there a simple, pedestrian
arithmetic sentence which is independent of
ZFC?'' preprint IEA--USP, to appear in {\it Synth\`ese} 
(2000). 
\bibitem {Johnson} D.\ S.\ Johnson, ``A catalogue of
complexity classes,'' in J.\ van Leeuwen, {\it Handbook of
Theoretical Computer Science}, Elsevier (1990). 
\bibitem {Sol} J.\ Ketonen, R.\ Solovay, ``Rapidly rising
Ramsey functions,'' {\it Annals of Mathematics} {\bf 113},
267 (1981). 
\bibitem {IM} S.\ C.\ Kleene, {\it Introduction to
Metamathematics}, Van Nostrand (1952). 
\bibitem {ML} S.\ C.\ Kleene, {\it Mathematical Logic}, John
Wiley (1967). 
\bibitem {Krei} G.\ Kreisel, ``On the interpretation of
non--finitist proofs,'' I and II, {\it J.\ Symbol.\ Logic}
{\bf 16}, 241 (1951) and {\bf 17}, 43 (1952). 
\bibitem {MY} M.\ Machtey and P.\ Young, {\it An
Introduction to the General Theory of Algorithms},
North--Holland (1978). 
\bibitem {Rog} H.\ Rogers Jr., {\it Theory of Recursive
Functions and Effective Computability}, McGraw--Hill (1967). 
\bibitem {Spen} J.\ Spencer, ``Large numbers and unprovable
theorems,'' {\it Amer.\ Math.\ Monthly} {\bf 90}, 669
(1983). 
\bibitem {Wai} S.\ S.\ Wainer, ``A classification of the
ordinal recursive functions,'' {\it Arch.\ Math.\ Logik} {\bf
13}, 136 (1970). 
\bibitem {Wai1} S.\ S.\ Wainer, ``Accessible recursive
functions,'' preprint, Univ.\ of Leeds, Dept.\ of Mathematics
(1999). 
\end {thebibliography}

\end {document}